\documentclass[leqno,12pt]{article}  
\usepackage{amsmath}
\usepackage{amssymb}

\usepackage[german,english]{babel}
\usepackage{amsthm}

\title{Rigid analytic spaces with overconvergent structure sheaf}

\author{\textsc{Elmar Grosse-Kl\"onne}}
\date{}

\theoremstyle{plain} 
\newtheorem{satz}{Theorem}[section]  
\newtheorem{lem}[satz]{Lemma}  
\newtheorem{pro}[satz]{Proposition}  
\newcommand{\ho}{\mbox{\rm Hom}}  
\newcommand{\ext}{\mbox{\rm Ext}} 
\newcommand{\spec}{\mbox{\rm Spec}}  

\newcommand{\spm}{\mbox{\rm Sp}}  
\newcommand{\spf}{\mbox{\rm Spf}}  

\newcommand{\bi}{\mbox{\rm im}}  
\newcommand{\ke}{\mbox{\rm Ker}}  
\newcommand{\koke}{\mbox{\rm Coker}}  

\newcommand{\q}{\mbox{\rm Frac}}  
\newcommand{\kara}{\mbox{\rm char}}  

\newcommand{\res}{\mbox{\rm res}}

\theoremstyle{remark}

\theoremstyle{definition}

\begin{document}
\maketitle
\footnote[0]{}
\footnote[0]{I would like to thank P.Schneider who suggested this topic for my thesis and showed constant interest in it. Besides, I thank him and A.Bertapelle for several useful conversations, and S.Bosch for pointing out an earlier inaccuracy of mine.}

\begin{abstract}

We introduce a category of 'rigid spaces with overconvergent structure sheaf' which we call dagger spaces --- this is the correct category in which de Rham cohomology in rigid analysis should be studied. 
We compare it with the (usual) category of rigid spaces, give Serre and Poincar\'{e} duality theorems and explain the relation with Berthelot's rigid cohomology.

\end{abstract}

%


\begin{center} {\bf Introduction} 
\end{center}

In rigid analysis the notion of overconvergent sheaf comes up in various shapes and for various reasons. The requirement that sections may slightly extend over boundaries should force sheaf properties usually shared only by sheaves on spaces whose topologies are given by 'wide open' subsets, or cohomological properties usually shared only by sheaf cohomology on 'spaces without boundary'. In this paper we establish a theory of 'rigid spaces with overconvergent structure sheaf' which we call dagger spaces. Continuing our reflection, the need for an overconvergent structure sheaf can be seen from the following examples: For a smooth rigid space $X$ which admits a closed immersion into a polydisc without boundary, the de Rham cohomology is (at least in many cases) finite dimensional, and there is a Serre duality formula, i.e. an interpretation of $\ext_{{\cal O}_X}^i({\cal O}_X,\omega_X)^{\bf \check{}}$ as some cohomology with compact support. But for a smooth rigid space which is affinoid, i.e. admits a closed immersion into a polydisc with boundary, the corresponding facts (with respect to its 'usual' structure sheaf) do not hold. Our point is that with an overconvergent structure sheaf things change.\\Let $k$ be a non-archimedean field. In the first section $k$-dagger algebras are introduced which form our substitutes for the $k$-affinoid algebras from (classical) rigid analysis (for which we refer to \cite{bgr}). 
A $k$-dagger algebra is a quotient of some algebra $W_n$ which is defined to be the subalgebra of the Tate algebra $T_n$ consisting of the power series with radius of convergence strictly greater than 1. 
On a $k$-dagger algebra $A$ there is a natural equivalence class of norms, and the completion $A'$ is a $k$-affinoid algebra. 
The functor $A\mapsto A'$ is studied in detail; for example, we find that the natural map $\tau:\spm(A')\to\spm(A)$ between the sets of maximal ideals is bijective. In section 2 we define dagger spaces. 
As in the rigid case one has a notion of affinoid subdomains of $\spm(A)$, and via $\tau$ these form a basis for the strong $G$-topology on $\spm(A')$. Imposing this $G$-topology on $\spm(A)$ one gets a locally $G$-ringed space, an affinoid $k$-dagger space. 
Then $k$-dagger spaces are built from affinoid ones. We define a faithful functor from the category of dagger spaces to the category of rigid spaces associating a rigid space $X'$ with a dagger space $X$ which has the same underlying $G$-topological space and the same stalks of structure sheaf. 
A smooth rigid space $Y$ admits an admissible open covering $Y=\cup V_i$ such that $V_i=U_i'$ for uniquely determined (up to isomorphism) dagger spaces $U_i$. 
Furthermore, this functor induces an equivalence between the respective subcategories formed by partially proper spaces (or: spaces without boundary). The third section is concerned with the cohomology of coherent modules. 
We prove the coherence theorem for proper morphisms assuming that $k$ is discretely valued. 
The fourth section deals with Serre duality for smooth affinoid dagger spaces, with Poincar\'{e} duality for the de Rham cohomology of smooth affinoid dagger spaces and of smooth (dagger or rigid) Stein spaces, and with the K\"unneth formula. 
In section 5 we compare our concept of overconvergence with the one developed by Berthelot (\cite{berco}), in particular we interpret his rigid cohomology as de Rham cohomology of certain dagger spaces.\\
In \cite{doktor},\cite{en2dag} we obtain, in case $k$ is discretely valued, for a big class of smooth dagger spaces, including the quasicompact ones, finiteness of de Rham cohomology, implying in particular finiteness of Berthelot's rigid cohomology.\\ 
Let us finally remark that, if $k$ is discretely valued, dagger spaces can be thought of as generic fibres of the weak formal schemes introduced in \cite{mer}, in the same way as rigid spaces can be thought of as generic fibres of admissible formal schemes (\cite{bolu}).\\

\section{Dagger algebras}
Let $k$ be a field complete with respect to a non-archimedean valuation $|.|$. We denote by $k_a$ its algebraic closure with value group $\Gamma^*=|k_a^*|=|k^*|\otimes{\mathbb Q}$.\\

\addtocounter{satz}{1}{\bf \arabic{section}.\arabic{satz}} \newcounter{teenro1}\newcounter{teenro2}\setcounter{teenro1}{\value{section}}\setcounter{teenro2}{\value{satz}} For $\rho\in{\mathbb R}^+$ define $$T_n(\rho) = \{\sum a_{\nu}X^{\nu}\in
k[[X_1,\ldots,X_n]]\quad | \quad |a_{\nu}|\rho^{|\nu|}\stackrel{|\nu|\to\infty}{\longrightarrow}0\},$$
where $|\nu|=\sum_{i=1}^n\nu_i$ for $\nu=(\nu_1,\ldots,\nu_n)\in{\mathbb N}^n$; in particular $T_n=T_n(1)$.\\
Every $T_n(\rho)$ is a $k$-affinoid algebra if $\rho\in\Gamma^*$, and the map $|.|_{\rho}:T_n(\rho)\to{\mathbb R}, \sum a_{\nu}X^{\nu}\mapsto \max|a_{\nu}|\rho^{|\nu|}$ defines a Banach norm on it (\cite{bgr},6.1.5).\\
Sometimes we write $k<\rho^{-1}.X_1,\ldots,\rho^{-1}.X_n>$ instead of $T_n(\rho)$ to specify names for the free variables.\\Also, if $A$ is a $k$-affinoid algebra, $f\in A$ and $\alpha\in\Gamma^*$, we write $A<\alpha^{-1}.f>$ instead of $A<a^{-1}.f^m>$, where $a\in k$ satisfies $|a|=\alpha^m, m\in{\mathbb N}$. This implies $\spm(A<\alpha^{-1}.f>)=\{x\in \spm(A)|\quad |f(x)|\le\alpha\}$.\\

\addtocounter{satz}{1}{\bf \arabic{section}.\arabic{satz}}  Following \cite{gue} we define the Washnitzer algebra $W_n$ to be$$W_n=k<X_1,\ldots,X_n>^{\dagger}=\cup_{\rho>1}T_n(\rho)=\cup_{\stackrel{\rho>1}{\rho\in\Gamma^*}}T_n(\rho).$$ 
We view $W_n$ as a subalgebra of $T_n$, equipped with the Gauss norm induced by the Gauss norm on $T_n$. A $k$-algebra $A$ is called a $k$-dagger algebra (or simply a dagger algebra if it is clear which is the ground field referred to) if there exists an $n\in{\mathbb N}$, an ideal $I<W_n$ and an isomorphism $A\cong W_n/I$ of $k$-algebras. A morphism of $k$-dagger algebras is a morphism of $k$-algebras. \\
By $\spm(A)$ we denote the set of maximal ideals of a dagger algebra $A$.\\

\addtocounter{satz}{1}{\bf \arabic{section}.\arabic{satz}} \newcounter{weier1} \setcounter{weier1}{\value{section}}\newcounter{weier2} \setcounter{weier2}{\value{satz}} A power series $g=\sum_{m=0}^{\infty}g_{m}(Y_1,\ldots,Y_{n-1})Y_n^{m}\in W_n=k<Y_1,\ldots,Y_n>^{\dagger}$ 
is called $Y_n$-distinguished of degree $k$ if $g_k$ is a unit in $W_{n-1}=k<Y_1,\ldots,Y_{n-1}>^{\dagger}$ and $|g|=|g_k|>|g_{m}|$ for all $m>k$.\\
A Weierstrass polynomial is a monic $\omega\in W_{n-1}[Y_n]\subset W_n$ with $|\omega|=1$.\\ 
The following (i), (ii) and (iv) are proven in \cite{gue}, whereas (iii) is a classical consequence of (i) (compare \cite{bgr},5.2.3/3,4):\\
(i) Let $g\in W_n$ be $Y_n$-distinguished of degree $k$. For all $f\in W_n$ there exist unique elements $q\in W_n, r\in W_{n-1}[Y_n]$ with $\deg(r)<k$ and $f=gq+r$.\\
(ii) Let $g\in W_n$ be $Y_n$-distinguished of degree $k$. There exists a unique Weierstrass polynomial $\omega\in W_{n-1}[Y_n]$ and a unit $e\in W_n$ with $g=e\omega$.\\
(iii) Let $\phi:W_n\to B$ be a finite morphisms of dagger algebras and suppose there exists a Weierstrass polynomial $\omega\in W_{n-1}[Y_n]\cap\ke(\phi)$. Then $W_{n-1}\to B$ is finite.\\
(iv) For all $0\ne f\in W_n$ there exists a $k$-algebra automorphism $\sigma$ of $W_n$ such that $\sigma(f)$ is $Y_n$-distinguished.\\

\addtocounter{satz}{1}{\bf \arabic{section}.\arabic{satz}}  
\newcounter{endrest1}\newcounter{endrest2}\setcounter{endrest1}{\value{section}}\setcounter{endrest2}{\value{satz}} As in \cite{bgr},5.2.6, 5.2.7 and 6.1.2, one deduces from \arabic{weier1}.\arabic{weier2}:\\(1) $W_n$ is a noetherian factorial Jacobson ring.\\
(2) Every dagger algebra $A$ admits a finite injection $W_n\to A$.\\
(3) If the nilradical of the ideal $q<A$ is a maximal ideal, then $k\hookrightarrow A/q$ is finite.\\
(4) Every ideal of $W_n$ is strictly closed in $W_n$, in particular closed.\\

\begin{pro}\label{spbij} $W_n\hookrightarrow T_n$ is faithfully flat. For every
maximal ideal $m<T_n$ the map $W_n/(m\cap W_n)\to T_n/m$ is bijective, in particular $m\cap W_n$ is maximal in $W_n$. For every
maximal ideal $y<W_n$ there are polynomials $p_1,\ldots,p_n\in
k[Y_1,\ldots,Y_n]$, each $p_i$ monic with respect to $Y_i$, which generate $y$. In particular, $W_n$ is regular, equidimensional of dimension $n$.\\
\end{pro}  

{\sc proof:} As sets one has $\spm(T_n)=\cap\spm(T_n(\rho))=\lim_{\leftarrow}\spm(T_n(\rho))$, for $\rho>1, \rho\to1$. On the other hand $\spec(W_n)=\lim_{\leftarrow}\spec(T_n(\rho))$ even as a projective limit of topological spaces for the Zariski topologies (\cite{EGA},IV,5.13.3, 8.2.9), 
and maximal ideals in $W_n$ induce inductive systems of maximal ideals in $(T_n(\rho))_{\rho}$ (use \arabic{endrest1}.\arabic{endrest2}). Together  surjectivity of $\spm(T_n)\to\spm(W_n)$ is implied. 
The claim on generation by polynomials now follows from the corresponding fact for $T_n$ (\cite{bgr},7.1.1/3).\\
To show flatness of $W_n\to T_n$ it is enough to show flatness of all localizations $(W_n)_m\to(T_n)_m$ in maximal ideals
$m$. But this follows from the bijectivity of the map of completed local rings (it is surjective because the associated map of graded rings is; it is injective for reasons of dimension).\\

\begin{pro}\label{normvgl}\label{stet} (1) For an ideal $I<W_n$ equip $W_n/I$ with the quotient semi-norm. Then this is a norm, and $\alpha:W_n/I\hookrightarrow T_n/I.T_n$ is the associated completion.\\
(2) Every $k$-algebra morphism $\phi:A=W_n/I\to B=W_m/J$ is continuous with respect to the topologies from (1).
\end{pro} 

{\sc proof:} Here all $W_n$ and all $T_n$ carry the Gauss norm.\\
(1) The first claim follows from \arabic{endrest1}.\arabic{endrest2}, the second from the fact that $W_n\hookrightarrow T_n$ is an isometry with dense image.\\
(2) We may suppose $A=W_n=k<X_1,\ldots,X_n>^{\dagger}$.\\ First we show that all elements $\phi(X_i)_{1\le i\le n}$ are power-bounded in $T_m/J.T_m$. By \cite{bgr},6.2.3/1, this is equivalent with $|\phi(X_i)|_{sup}\le1$ in $T_m/J.T_m$ for all $1\le i\le n$. 
But this is equivalent with $|\phi(X_i)|_{sup}\le1$ in $B$ for all $1\le i\le n$, because $B\to T_m/J.T_m$ induces a bijection between the sets of maximal ideals and isomorphisms of residue fields in maximal ideals. For a maximal ideal
$m<B$ and $\nu\ge1$ consider the canonical map $$\psi:W_n\to B\to
B/(m^{\nu})$$ and the induced map $\bar{\psi}:W_n/\ke(\psi)\hookrightarrow B/(m^{\nu})$. By \arabic{endrest1}.\arabic{endrest2} the respective quotient semi-norms on $W_n/\ke(\psi)$ and $B/(m^{\nu})$ are norms, and since $B/(m^{\nu})$ is finite dimensional (by \arabic{endrest1}.\arabic{endrest2}) so is $W_n/\ke(\psi)$, 
so $\bar{\psi}$ is continuous for these quotient (semi-)norms (use \cite{bgr},2.3.3). Consequently, $\psi$ is continuous, implying $$(*)\quad\phi(g(X))\equiv g(\phi(X)\quad\mod\quad
m^{\nu})$$ for all $g\in W_n$ and in particular (with $\nu=1$) the convergence of
$g(\phi(X) \mod m)$ in $B/m$. If $|\phi(X_i)|_{B/m}=\delta>1$, then $g(\phi(X_i) \mod m)$ would not converge for $g\in W_n-T_n(\delta)$. It follows $|\phi(X_i)|_{B/m}\le 1$ for all $i$.\\
Now since all elements $\phi(X_i)_{1\le i\le n}$ are power-bounded in $T_m/J.T_m$, one can define a continuous map $\phi':T_n\to T_m/J.T_m$ by requiring $\phi'(X_i)=\phi(X_i)$. 
Then from $(*)$ and the identity $\phi'(g(X))=g(\phi'(X))$ for all $g\in T_n$ one gets: For every $g\in W_n$ the congruences $\phi(g(X))\equiv
\phi'(g(X)) \mod (m^{\nu})$ of $T_m/J.T_m$-elements hold (all $\nu\ge1$, all maximal ideals $m<B$). This means $\phi(g(X))=\phi'(g(X))$ because the intersection of all ideals $m^{\nu}$ is the zeroideal (\cite{bgr},6.1.3).\\

In the sequel we provide every dagger algebra with a norm of the equivalence class described in \ref{normvgl}.\\

\begin{satz}\label{eigueb} Let $A$ be a dagger algebra with completion $\tau:A\to A'$. \\(1) $\tau$ ist faithfully flat and induces a bijection between the sets of maximal ideals.\\
(2) $\tau$ induces isomorphisms between the completions in maximal ideals.\\
(3) $A$ is reduced (resp. normal, resp. regular) if and only if $A'$ is.\\(4) If $A$ is reduced, the supremum semi-norm on $A$ is a norm and belongs to the equivalence class described in \ref{normvgl}.
\end{satz}

{\sc Proof:} (1) follows from \ref{spbij}.\\(2): Write $A=W_n/I$. If $m<W_n$ is a maximal
ideal contatining $I$ and $m'<T_n$ is the corresponding maximal ideal containing $I.T_n$, one finds $(W_n/I)\hat{}_m=(W_n)\hat{}_m\otimes_{W_n}(W_n/I)=(T_n)\hat{}_{m'}\otimes_{W_n}(W_n/I)=(T_n)\hat{}_{m'}\otimes_{T_n}(T_n/I.T_n)=(T_n/I.T_n)\hat{}_{m'}$ for the completed local rings (compare the proof of \ref{spbij}).\\
(3) follows from (2), because these properties can be checked in the completed local rings (all rings are excellent; apply \cite{EGA},IV, 7.8).\\(4): By (1) and (2) $\alpha$ is an isometry with respect to the supremum semi-norms, so it is enough to quote the analogous statement for $A'$ (which is reduced by (3)).\\ 

\begin{lem} \label{extens} Let $\psi:A\to B$ be a morphism of dagger algebras. Suppose we are given representations $A=W_n/(f_1,\ldots,f_r)$ and $B=W_m/(g_1,\ldots,g_s)$ with $f_1,\ldots,f_r\in T_n(\lambda_A)$ and $g_1,\ldots,g_s\in T_m(\lambda_B)$ for $\lambda_A,\lambda_B>1$ in $\Gamma^*$. 
For any $\lambda_A\ge\lambda>1$ and $\lambda_B\ge\lambda'>1$ in $\Gamma^*$ we put $A_{\lambda}=T_n(\lambda)/(f_1,\ldots,f_r)$ and $B_{\lambda'}=T_m(\lambda')/(g_1,\ldots,g_s)$. 
Then for any such $\lambda$ there is a $\lambda'$ such that the composition $A_{\lambda}\to A\to B$ factorizes over $B_{\lambda'}\to B$.
\end{lem}

{\sc Proof:} Let $c\in k^*, s\in {\mathbb N}$ with $\lambda^s=|c^{-1}|$. Then $|c.X_i^l|_{\lambda}\le1$ for all $i=1,\ldots,n$, all $l=0,\ldots,s$, where the norm $|.|_{\lambda}$ 
on $T_n(\lambda)=k<\lambda^{-1}.X_1,\ldots,\lambda^{-1}.X_n>$ is defined as in \arabic{teenro1}.\arabic{teenro2}. 
The continuous map $$T_{(s+1)n}=k<T_{il}>_{\stackrel{i=1,\ldots,n}{l=0,\ldots,s}}\stackrel{p}{\to} T_n(\lambda), \quad T_{il}\mapsto c.X_i^l$$ is surjective (see the proof of \cite{bgr},6.1.5/4). 
Denote by $\phi_{\lambda}:A_{\lambda}\to A$ and $\beta_{\lambda}:T_n(\lambda)\to A_{\lambda}$ the canonical maps. The elements $\phi_{\lambda}(\beta_{\lambda}(p(T_{il})))$ are topologically nilpotent, and so are the elements $\psi(\phi_{\lambda}(\beta_{\lambda}(p(T_{il}))))\in B$. 
Therefore there exist $\lambda''>1$ and power-bounded elements $t_{il}$ ($i=1,\ldots,n$ and $l=0,\ldots,s$) in $B_{\lambda''}$ with $\gamma_{\lambda''}(t_{il})=\psi(\phi_{\lambda}(\beta_{\lambda}(p(T_{il}))))$, where 
$\gamma_{\lambda''}:B_{\lambda''}\to B$ is the canonical map. For the continuous map$$q:T_{(s+1)n}\to B_{\lambda''},\quad T_{il}\mapsto t_{il}$$
one has $\gamma_{\lambda''}\circ q=\psi\circ\phi_{\lambda}\circ\beta_{\lambda}\circ p$, hence $\gamma_{\lambda''}(q(\ke(\beta_{\lambda}\circ p)))=0$ 
and therefore$$\gamma_{\lambda''}(q(\ke(\beta_{\lambda}\circ p)).B_{\lambda''})=(0)<B.$$
Let $q(\ke(\beta_{\lambda}\circ p)).B_{\lambda''}=(r_1,\ldots,r_m)<B_{\lambda''}$. Since $B$ is the direct limit of the maps $\mu_{\eta'\eta}:B_{\eta'}\to B_{\eta}$ ($\eta'\ge\eta>1$), 
there exists a $\lambda'$ with $\lambda''\ge\lambda'>1$ such that $\mu_{\lambda''\lambda'}(r_j)=0$, all $j=1,\ldots,m$ and we get the wanted factorization $A_\lambda\to B_{\lambda'}$.\\

\addtocounter{satz}{1}{\bf \arabic{section}.\arabic{satz}}  We call a normed $k$-algebra $A$ weakly complete (relative $k$), if for every collection $x_1,\ldots,x_n\in A$ of power-bounded elements the $k$-algebra morphism $k[X_1,\ldots,X_n]\to A, X_i\mapsto x_i$ 
admits a continuous extension $k<X_1,\ldots,X_n>^{\dagger}\to A$.\\
$k$-dagger algebras are weakly complete: If $A=W_k/I$ and $x_1,\ldots, x_n\in A$ are power-bounded, let $T_n=k<X_1,\ldots,X_n>\stackrel{\psi}{\to}T_k/I.T_k$ be the continuous map sending $X_i\mapsto x_i$. 
We have to show that for all $\lambda>1$ the composition $T_n(\lambda)\to T_n\stackrel{\psi}{\to}T_k/I.T_k$ factorizes over $T_k(\lambda')/I'\to T_k/I.T_k$ for appropriate 
$\lambda'=\lambda'(\lambda)>1$ and $I'<T_k(\lambda')$ with $I'.W_k=I$. This can be seen as in the proof of \ref{extens}.\\ 

\begin{lem}\label{enddag} Let $A$ be a dagger algebra and $A\to B$ a finite $k$-algebra morphism. Then $B$ is a dagger algebra.
\end{lem}  

{\sc Proof:} We may suppose $A=W_n/I.W_n$ and $B=(W_n/J.W_n)[X]/(f)$ with ideals $I\subset J\subset T_n(\lambda_0)$ and a monic polynomial $f\in T_n(\lambda_0)[X]$ for some $\lambda_0>1$. 
Since $T_n(\lambda_0)/J.T_n(\lambda_0)\to T_n(\lambda_0)[X]/(J+f)$ is finite we may (after changing $f$ if necessary) also suppose that $X$ induces a power-bounded element in $T_n(\lambda_0)[X]/(J+f)$. 
Then
$$D_{\lambda}=T_n(\lambda)[X]/(J+f)\to T_n(\lambda)<\lambda^{-1}.X>/(J+f)=T_{n+1}(\lambda)/(J+f)=E_{\lambda}$$is bijective for every $1\le\lambda\le\lambda_0$ 
and hence also$$B=\lim_{\stackrel{\to}{\lambda>1}}D_{\lambda}\to\lim_{\stackrel{\to}{\lambda>1}}E_{\lambda}=W_{n+1}/(J+f).$$

\addtocounter{satz}{1}{\bf \arabic{section}.\arabic{satz}}\newcounter{modnor1}\newcounter{modnor2}\setcounter{modnor1}{\value{section}}\setcounter{modnor2}{\value{satz}} For every dagger algebra $A$ the finite $A$-modules can be canonically equipped with an equivalence class of norms such that all morphisms of finite $A$-modules become continuous. 
Namely, if $M$ is a finite $A$-module and if $A'$ is the completion of $A$, take the restriction of the complete $A'$-module norm on $M\otimes_AA'$ from \cite{bgr},3.7.3/3, to $M$. The claim on continuity then follows from \cite{bgr},3.7.3/2.\\  

\begin{satz}\label{hauptend} Let $\phi:A\to B$ be a morphism of dagger algebras and $\phi':A'\to B'$ its completion.\\(a) $\phi$ is surjective (resp. bijective) if and only $\phi'$ is.\\
(b) The following are equivalent:\\(i) $\phi$ is finite.\\(ii) $\phi'$ is finite.\\(iii) $\tilde{\phi}:\tilde{A}\to\tilde{B}$ is finite.\\(c) If $\phi$ finite, then $\phi$ is injective if and only if $\phi'$ is.\\
(d)\label{normend} Let $\phi$ be finite. On $B$ the equivalence classes of norms from \arabic{modnor1}.\arabic{modnor2} and from \ref{normvgl} coincide, and the canonical map $A'\otimes_AB\to B'$ is an isomorphism.\\(e) If $\phi$ is finite, $\phi$ and $\phi'$ are strict.\\
\end{satz}

{\sc Proof:} Here we use the usual notation $\tilde{R}=R^0/R^{00}$ for a normed ring $R$, see \cite{bgr}.\\
(a) and (c) follow from (b) and (d) ($A\to A'$ is faithfully flat).\\
(d): $A'\otimes_AB$ is a $k$-affinoid algebra with $k$-Banach algebra norm equivalent to the norm it carries as a finite $A'$-module (\cite{bgr},3.7.3/3,
3.7.4/1). From (b) one gets that $A'\otimes_AB\to B'$ is finite and by \cite{bgr},3.7.4/1, this means that the $k$-Banach algebra norm on $B'$ is a $A'\otimes_AB$-algebra norm. This gives the first claim, and then also the second because $A'\otimes_AB$ is the completion of $B$ with respect to its norm as a finite $A$-module.\\
(e): $\phi'$ is finite by (b), so we can apply \cite{bgr},1.1.9/2, 6.1.3/4.\\  
(b): Setting $\lambda=1$ in the proof of \ref{enddag} one obtains (b),(i)$\Rightarrow$(b),(ii), and (b),(ii)$\Rightarrow$(b),(iii) is \cite{bgr},6.3.4/2. Finally, (b),(iii)$\Rightarrow$(b),(i) can be proven literally as in \cite{bgr},6.3.2/1,2.\\

\begin{pro}\label{artapp} (\cite{boartin}) Let $I<W_n$ be an ideal, put $A=W_n/I\subset A'=T_n/(I.T_n)$ and fix a Banach norm $|.|_{A'}$ on $A'$. Let $(f_i)_{i\in I}\in A<Y_1,\ldots,Y_m>^{\dagger}=W_{n+m}/I.W_{n+m}$ 
and suppose there exist $\bar{y}_1\ldots,\bar{y}_m\in A'$ with $|\bar{y}_j|_{sup}\le 1$ and $f_i(\bar{y})=0$ for all $i\in I$. Then for all $\delta>0$ there exist $y_1,\ldots,y_m\in A$ with $|y_j|_{sup}\le 1$ and $f_i(y)=0$ for all $i\in I$, satisfying $|y-\bar{y}|_{A'}<\delta$.  
\end{pro}  

\begin{pro}\label{nixend} Let $A$ be a dagger algebra, $A\to A'$ its completion.\\(1) If $A$ is reduced, then it is integrally closed in $A'$.\\(2) $A$ is an integral domain if and only $A'$ is.
\end{pro} 

{\sc Proof:} (1) follows from \ref{artapp}, see \cite{boartin},sect.2.\\(2) Let $A$ be an integral domain. As in \cite{bgr},6.1.2/4, one sees that $A$ is japanese, so the normalization $A\to B$ is finite. 
By \ref{enddag} and \ref{hauptend} this means that also the completion $A'\to B'$ is injective, so it is enough to show that $B'$ is an integral domain. 
Since it is normal, all its local rings are integral domains. On the other hand, the connectedness of $\spm(B)$ (for the Zariski topology) implies the connectednes of $\spm(B')$ (for the Zariski topology, or equivalent: for the strong $G$-topology; use \cite{fies},3.3.3).\\ 

\begin{lem} \label{abgimapp} \label{einddag} Let $p_1:A_1\to B_1$ and $p_2:A_2\to B_2$ be surjective morphisms of dagger algebras with completions $p_1':A_1'\to B_1'$ and $p_2':A_2'\to B_2'$. 
Let furthermore $\phi_A:A_1'\to A_2'$ and $\phi_B:B_1'\to B_2'$ be morphisms satisfying $\phi_B\circ p_1'=p_2'\circ\phi_A$. 
Then there exist morphisms $\gamma_A:A_1\to A_2$ and $\gamma_B:B_1\to B_2$ satisfying $\gamma_B\circ p_1=p_2\circ\gamma_A$. 
If $\epsilon>0$ and a Banach norm $|.|$ on $A'_2$ are given, $\gamma_A$ can be chosen such that for the completion $\gamma_A'$ one has $|\gamma_A'-\phi_A|<\epsilon$. 
In particular, if $\phi_A$ is an isomorphims, $\gamma_A$ can be chosen to be an isomorphism.   
\end{lem}

{\sc Proof:} Write $A_1=k<X_1,\ldots,X_n>^{\dagger}/(f_1,\ldots,f_r)$ and $B_1=k<X_1,\ldots,X_n>^{\dagger}/(f_1,\ldots,f_r,g_1,\ldots,g_l)$ (according to the surjection $p_1$) and $\ker(p_2)=(d_1,\ldots,d_p)$ with $d_i\in A_2$. Consider solutions $(x_i,y_{jk})_{i,j,k}\in(A_2')^{n+lp}$ 
(where $i\in\{1,\ldots,n\}, j\in\{1,\ldots,l\}$ and $k\in\{1,\ldots,p\}$) of the system$$f_q(x_1,\ldots,x_n)=0\mbox{ for }q\in\{1,\ldots,r\}$$
$$g_j(x_1,\ldots,x_n)=\sum_{k=1}^py_{jk}.d_k\mbox{ for }j\in\{1,\ldots,l\}.$$The pair $(\phi_A,\phi_B)$ corresponds to such a solution. 
By \ref{artapp} it can be approximated by a solution in $A_2^{n+lp}$, defining $(\gamma_A,\gamma_B)$. If $\phi_A$ is an isomorphim so is $\gamma_A'$ (if $|\gamma_A'-\phi_A|<1$). Now apply \ref{hauptend}.\\

\addtocounter{satz}{1}{\bf \arabic{section}.\arabic{satz}} \newcounter{tensor1}\newcounter{tensor2}\setcounter{tensor1}{\value{section}}\setcounter{tensor2}{\value{satz}} In the category of dagger algebras tensor products exist. In fact, suppose we are given surjections $f_1:W_{n_1}=k<X_1,\ldots,X_{n_1}>^{\dagger}\to A_1$ and $f_2:W_{n_2}=k<Y_1,\ldots,Y_{n_2}>^{\dagger}\to A_2$ and moreover morphisms $\alpha_1:B\to
A_1$, $\alpha_2:B\to A_2$ of dagger algebras. Define $A_1\otimes_B^{\dagger}A_2$ to be the image of $W_{n_1+n_2}$ under the map $\psi:T_{n_1+n_2}\to
A_1\hat{\otimes}_BA_2=A_1'\hat{\otimes}_{B'}A_2'$ into the completed tensor product of the $k$-affinoid algebras $A_1'$ and $A_2'$ over the $k$-affinoid algebra $B'$, the completion of $B$. The universal property of  $A_1\otimes_B^{\dagger}A_2$ is then deduced from the universal property of $A_1\hat{\otimes}_BA_2=A_1'\hat{\otimes}_{B'}A_2'$ together with the weak completeness of dagger algebras.\\

\addtocounter{satz}{1}{\bf \arabic{section}.\arabic{satz}} For $\lambda\in\Gamma^*$ we write$$k<\lambda^{-1}.Z_1,\ldots,\lambda^{-1}.Z_n>^{\dagger}=\cup_{\rho>\lambda}k<\rho^{-1}.Z_1,\ldots,\rho^{-1}.Z_n>.$$
This is a $k$-dagger algebra. If $A$ is a dagger algebra and $Z$ is a free variable, we write $A<\lambda^{-1}.Z>^{\dagger}=A\otimes_k^{\dagger}k<\lambda^{-1}.Z>^{\dagger}$ and in particular $A<Z>^{\dagger}=A<1^{-1}.Z>^{\dagger}$. 
Note that in general $A<\lambda^{-1}.Z>^{\dagger}$ is strictly smaller then the algebra of all power series $f=\sum
a_{\nu}Z^{\nu}\in A[[Z]]$ with $|a_{\nu}|\rho^{\nu}\rightarrow0$ for some $\rho=\rho(f)>\lambda$. Note also that the notation is consistent with finite extensions $k\subset k'$ of the ground field.\\

\addtocounter{satz}{1}{\bf \arabic{section}.\arabic{satz}} \newcounter{ratio1}\newcounter{ratio2}\setcounter{ratio1}{\value{section}}\setcounter{ratio2}{\value{satz}} Let $A$ be a dagger algebra and suppose the elements $g,f_1,\ldots,f_m\in A$ have no common zero. Set$$A<f/g>^{\dagger}=A<X_i>_{1\le i\le m}^{\dagger}/((gX_i-f_i)_{1\le i\le
m}).$$ Then, in the category of dagger algebras, $A\to A<f/g>^{\dagger}$ satisfies the universal property analogous to \cite{bgr},6.1.4.\\

\section{Dagger spaces}
\addtocounter{satz}{1}{\bf \arabic{section}.\arabic{satz}}  Let $A$ be a dagger algebra, $X=\spm(A)$. A subset $U\subset X$ is called an affinoid subdomain of $X$, if there exists a dagger algebra $B={\cal O}_X(U)$ and a morphism of dagger algebras 
$\pi:A\to B$ with $\bi(\spm(\pi))\subset U$ such that for all dagger algebras $D$ the map $$\ho_k(B,D)\to\{f\in\ho_k(A,D)| \mbox{ }\bi(\spm(f))\subset U\},\quad g\mapsto g\circ\pi$$ is bijective. 
(Then it follows that $B$ is uniquely determined and that $\spm(B)=U$). Every $A<f/g>^{\dagger}$ as in \arabic{ratio1}.\arabic{ratio2} defines an affinoid subdomain: these are called rational subdomains (and if $g=1$ they are called Weierstrass domains). These properties are transitive and stable under arbitrary base change $A\to D$.\\

\addtocounter{satz}{1}{\bf \arabic{section}.\arabic{satz}} If the morphism of dagger algebras $A\to B$ defines an affinoid subdomain of $\spm(A)$, then the completion $A'\to B'$ defines an affinoid subdomain of the rigid affinoid space $\spm(A')$ (and the corresponding injections of sets of maximal ideals are naturally identified). 
To see this, prove the analogues of \cite{bgr},7.2.2/1, in the dagger context, compare the completed local rings and apply \cite{bgr},7.3.3/5, 8.2.1/4.\\

\addtocounter{satz}{1}{\bf \arabic{section}.\arabic{satz}}  For $x\in X=\spm(A)$ set ${\cal O}_{X,x}=\lim_{\to}{\cal O}_X(U)$, where $U$ runs through the set of all affinoid subdomains of $X$ with $x\in U$. 
We say that a morphism of dagger algebras $\pi:A\to B$ defines an open immersion (resp. a locally closed immersion) if for all $x\in\spm(B)$ the induced maps ${\cal O}_{\spm(A),y}\to{\cal O}_{\spm(B),x}$ are bijective 
(resp. surjective); here $y=\spm(\pi)(x)$.\\

\addtocounter{satz}{1}{\bf \arabic{section}.\arabic{satz}} Let $A$ be a dagger algebra, $A'$ its completion and $x$ an element of the canonically identified sets $X=\spm(A)$ and $X'=\spm(A')$. 
Then the canonical morphisms ${\cal O}_{X,x}\to {\cal O}_{X',x}$ are isomorphisms. 
In fact, to show surjectivity write $A=W_n/I.W_n$ and $X'(\rho)=\spm(T_n(\rho)/I)$ for some $\rho>1$ and ideal $I<T_n(\rho)$. Then one has ${\cal O}_{X'(\rho),x}\cong{\cal O}_{X',x}$, 
and for $g\in{\cal O}_{X'(\rho),x}$ there exist $f_1,\ldots,f_r\in x<T_n(\rho)/I$ and an element $\dot{g}\in (T_n(\rho)/I)<f_1,\ldots,f_r>$ which induces $g$. But $\dot{g}$ induces also an element of $(W_n/I.W_n)<\rho.f_1,\ldots,\rho.f_r>^{\dagger}$ and hence the wanted element in ${\cal O}_{X,x}$.\\
We deduce that a morphism of dagger algebras defines an open (resp. locally closed) immersion if and only if the associated morphism of $k$-affinoid algebras does.

\begin{pro}\label{GG} Let $A$ be a dagger algebra, $A'$ its completion, $X=\spm(A)$ and $X'=\spm(A')$. 
If $A'\to D$ defines a locally closed immersion of rigid affinoid spaces, there exists a finite covering $X=\cup X_i$ by rational subdomains $X_i=\spm(A<f_i/g_i>^{\dagger})$ 
with the following property: If $A'<f_i/g_i>=D_i'$ are the completions, 
the morphisms $D_i'\to D_i'\hat{\otimes}_{A'}D$ define Runge immersions of affinoid rigid spaces.
\end{pro}

{\sc Proof:} If $A=W_m/I.W_m$ for some $\rho>1$ and ideal $I<T_m(\rho)$, then $E=T_m(\rho)/I\to D$ defines a locally closed immersion. By \cite{bgr},7.3.5/1, there exists a covering of $\spm(E)$ 
by finitely many rational subdomains $W_i$ such that all morphisms $\spm(D)\cap W_i\to W_i$ are Runge immersions. This induces the desired covering of $X$.\\ 

Note that a morphism of dagger algebras defines a Runge immersion (that is: factorizes into a Weierstrass domain and a closed immersion) if and only if this is true for the associated morphism of affinoid rigid spaces.\\

\begin{pro}\label{azyk}\label{modazyk} Let $A$ be a dagger algebra and $M$ an $A$-module. Every finite covering $X=\spm(A)=\cup_{i=1}^nU_i$ 
by affinoid subdomains (we call such coverings affinoid coverings) is acyclic for the presheaf $U\mapsto{\cal O}_X(U)\otimes_AM$.
\end{pro} 

{\sc Proof:} As in \cite{bgr},8.2.1/5, this is reduced to the case $M=A$. With \ref{GG} at hand one can then follow the proof of \cite{bgr},8.2.1/1.\\

\addtocounter{satz}{1}{\bf \arabic{section}.\arabic{satz}} As in \cite{bgr},8.2.1, we deduce from \ref{azyk} that if a morphism of dagger algebras $\pi:A\to B$ defines an open immersion, then $\bi(\spm(\pi))$ is an affinoid subdomain of $\spm(A)$. 
It follows that $\pi$ defines an affinoid subdomain if and only if this is true for the completion $\pi':A'\to B'$.\\
If $A$ is a dagger algebra, $A'$ its completion, we provide the set $\spm(A)$ with the $G$-topology which is induced from the ('strong') $G$-topology on $\spm(A')$ via the canonical bijection $\spm(A)\cong\spm(A')$.\\

\begin{pro}\label{gleichtop} Let $A$ be a dagger algebra. The affinoid subdomains of $\spm(A)$ form a basis for the $G$-topology. 
For any admissible covering of an affinoid subdomain of $\spm(A)$ there exists a finite refinement consisting of affinoid subdomains.
\end{pro}

{\sc Proof:} This follows from the corresponding statements for the completion $A'$ as soon as one has: Via the bijection $\spm(A)\cong\spm(A')$ every affinoid subdomain $U$ of $\spm(A')$ is  
the union of finitely many affinoid subdomains of $\spm(A)$. Using \ref{GG} one reduces the problem to the case where $U\subset\spm(A')$ is a Weierstrass domain and then concludes by the density of $A$ in $A'$.\\

\addtocounter{satz}{1}{\bf \arabic{section}.\arabic{satz}} One could instead have defined the $G$-topology on $\spm(A)$ in complete analogy with the rigid case as in \cite{bgr}, obtaining a category of dagger spaces as good as we obtain it here. 
The reason for our choice of $G$-topology (which is coarser than the just mentioned one) is that we want to have a functor from the category of dagger spaces to the category of rigid spaces.\\

\addtocounter{satz}{1}{\bf \arabic{section}.\arabic{satz}} Let $A$ be a dagger algebra, $X=\spm(A)$ and $M$ an $A$-module. Because of \ref{modazyk}, \ref{gleichtop} the presheaf $U\mapsto M\otimes_A{\cal O}_X(U)$ on affinoid subdomains $U\subset\spm(A)$ extends in a unique way to a sheaf  $M\otimes{\cal O}_X$ for the $G$-topology 
(use \cite{bgr},9.2.3/1). $A\otimes{\cal O}_X={\cal O}_X$ is naturally a sheaf of $k$-algebras, and all $M\otimes{\cal O}_X$ are modules over ${\cal O}_X$.\\

\addtocounter{satz}{1}{\bf \arabic{section}.\arabic{satz}}\newcounter{volltreu1}\newcounter{volltreu2}\setcounter{volltreu1}{\value{section}}\setcounter{volltreu2}{\value{satz}}  The assignement $A\mapsto(\spm(A)=X,{\cal O}_X)$ is a fully faithful functor from the category of $k$-dagger algebras into the category of locally $G$-ringed spaces over $k$. 
A morphism $A\to B$ of dagger algebras defines an open immersion if and only if $\spm(B)\to\spm(A)$ is an open immersion of $G$-topological spaces.\\

\addtocounter{satz}{1}{\bf \arabic{section}.\arabic{satz}}  A locally $G$-ringed space over $k$ in the essential image of the functor in \arabic{volltreu1}.\arabic{volltreu2} is called an affinoid $k$-dagger space. 
A locally $G$-ringed space $(X,{\cal O}_X)$ over $k$ is called a $k$-dagger space, if it satisfies the following conditions:\\
(i) The $G$-topology is saturated, that is (G0)(G1)(G2) from \cite{bgr},9.2.1, hold.\\(ii) There exists an admissible covering $X=\cup_iU_i$ such that all $(U_i,{\cal O}_X|_{U_i})$ are affinoid $k$-dagger spaces.\\
A morphism of dagger spaces is a morphism of locally $G$-ringed spaces.\\Again we often omit the $k$ and write of (affinoid) dagger spaces.\\

\addtocounter{satz}{1}{\bf \arabic{section}.\arabic{satz}}  Open subspaces of dagger spaces are dagger spaces. The category of dagger spaces has fibre products (use \arabic{tensor1}.\arabic{tensor2}). If $X$ is a $k$-dagger space and $k\subset k'$ is a finite field extension, then $X\times_{\spm(k)}\spm(k')$ is in a natural way a $k'$-dagger space.\\We denote by ${\bf D}=\{x\in k| |x|\le1\}$ (resp. ${\bf D}^0=\{x\in k| |x|<1\}$) the unit disc with (resp. without) boundary, with its dagger or rigid structure.\\

\addtocounter{satz}{1}{\bf \arabic{section}.\arabic{satz}}  Let $X$ be a dagger space. An ${\cal O}_X$-module ${\cal F}$ is called coherent, if there exists an affinoid covering $X=\cup_iU_i=\cup_i\spm(A_i)$, finite $A_i$-modules 
$M_i$ and isomorphisms ${\cal F}|_{U_i}\cong M_i\otimes{\cal O}_{U_i}$ of ${\cal O}_{U_i}$-modules.\\

\begin{lem}\label{mohomex} Let $\rho_0>1$, $I<T_n(\rho_0)$, $A=W_n/I.W_n$ and $B_{\rho}=T_n(\rho)/I.T_n(\rho)$ for $\rho_0\ge\rho>1$.\\(1) For any finite $A$-module $M$ there exists $\rho$, 
a finite $B_{\rho}$-module $M_{\rho}$ and an isomorphism $M\cong M_{\rho}\otimes_{B_{\rho}}A$.\\
(2) For any two finite $B_{\rho}$-modules $(M_1)_{\rho}, (M_2)_{\rho}$ and any $A$-linear morphism $\nu:(M_1)_{\rho}\otimes_{B_{\rho}}A\to(M_2)_{\rho}\otimes_{B_{\rho}}A$ there exists, after shrinking $\rho$ if necessary, a 
$B_{\rho}$-linear morphism $\nu_{\rho}:(M_1)_{\rho}\to(M_2)_{\rho}$ with $\nu=\nu_{\rho}\otimes_{B_{\rho}}A$. If $\nu$ is an isomorphism, $\nu_{\rho}$ can be chosen to be an isomorphism.\\
\end{lem}

{\sc Proof:} (1) Choose a finite presentation $A^r\stackrel{\psi}{\to}A^s\to M\to0$. Because of $A^s=\lim_{\stackrel{\to}{\rho>1}}B_{\rho}^s$ there is a $B_{\rho}$-linear map $\psi_{\rho}:B_{\rho}^r\to B_{\rho}^s$ such that $\psi_{\rho}\otimes_{B_{\rho}}A=\psi$. Take $M_{\rho}=\koke(\psi_{\rho})$.\\
(2) Choose finite presentations $B_{\rho}^{r_i}\stackrel{(\psi_i)_{\rho}}{\to}B_{\rho}^{s_i}\to(M_i)_{\rho}\to0$. Then $\nu$ is induced by $A$-linear maps 
$\phi_s:A^{s_1}\to A^{s_2}$, $\phi_r:A^{r_1}\to A^{r_2}$ which form a commutative diagram with the maps $(\psi_i)_{\rho}\otimes_{B_{\rho}}A$. Now conclude as in (1).\\

\begin{satz}\label{kiehl} Let $X=\spm(A)$ be an affinoid dagger space and ${\cal F}$ a coherent ${\cal
O}_X$-modules. Then there exists a finite $A$-module $M$ and an isomorphism ${\cal F}\cong M\otimes{\cal O}_X$ of ${\cal
O}_X$-modules.
\end{satz}  

{\sc Proof:} Using \ref{GG} we may suppose (compare the proof of \cite{bgr},8.4.3/3) that there exist $f_1, f_2\in A$, $f_2=f_1^{-1}$, finite $A_i$-modules $M_i$ and isomorphisms 
${\cal F}|_{U_i}\cong M_i\otimes{\cal O}_{U_i}$ of ${\cal O}_{U_i}$-modules, where we put $A_i=A<f_i>^{\dagger}$ and $U_i=\spm(A_i)$.\\
Write $A=W_n/I.W_n$ for some $I<T_n(\rho_0)$ and $\rho_0>1$, and suppose $f_1, f_2\in A$ are induced by $f_1, f_2\in T_n(\rho_0)/I$. Setting $A_{12}=A<f_1,f_2>^{\dagger}$, $A_{\rho}=T_n(\rho)/(I)$, $$A_{i,\rho}=T_n(\rho)<\rho^{-1}.Y>/(I+(Y-f_i))$$ 
and $$A_{12,\rho}=T_n(\rho)<\rho^{-1}.Y_1,\rho^{-1}.Y_2>/(I+(Y_1-f_1)+(Y_2+f_2))$$ 
for $\rho_0\ge\rho>1$ there exist by \ref{mohomex} finite $A_{i,\rho}$-modules $M_{i,\rho}$ (for some $\rho$) inducing the $M_i$. 
By the sheafproperty of ${\cal F}$ there exists an isomorphism $$(M_{1,\rho}\otimes_{A_{1,\rho}}A_{12,\rho})\otimes_{A_{12,\rho}}A_{12}\cong(M_{2,\rho}\otimes_{A_{2,\rho}}A_{12,\rho})\otimes_{A_{12,\rho}}A_{12}$$
of $A_{12}$-modules, which, after shrinking $\rho$ if necessary (\ref{mohomex}), is induced by an isomorphism $(M_{1,\rho}\otimes_{A_{1,\rho}}A_{12,\rho})\cong(M_{2,\rho}\otimes_{A_{2,\rho}}A_{12,\rho})$. 
Now apply \cite{bgr},9.4.3/3, to get a finite $A_{\rho}$-module $M_{\rho}$ such that $M_{i,\rho}\cong M_{\rho}\otimes_{A_{\rho}}A_{i,\rho}$, which means $M_i=M\otimes_AA_i$ for the finite $A$-module $M=M_{\rho}\otimes_{A_\rho}A$.\\

\addtocounter{satz}{1}{\bf \arabic{section}.\arabic{satz}}  Now all the analogous assertions of \cite{bgr},9.5, 9.6, 9.6.2, on coherent ideals, closed immersions, separated and finite morphisms (analogous definitions) can be literally translated to the dagger context.\\

\addtocounter{satz}{1}{\bf \arabic{section}.\arabic{satz}}  A  rigid space $X$ is called quasialgebraic, if there exists an admissible covering $X=\cup_i X_i$ and for all $i$ an open immersion $X_i\to Y_i^{\mbox{an}}$ into the analytification $Y_i^{\mbox{an}}$ of an affine $k$-scheme 
$Y_i$ of finite type. If $X$ is quasialgebraic, there exists an admissible affinoid covering $X=\cup_i X_i$, $n_i, r_i\in{\mathbb N}$, polynomials $f_{ij}\in k[X_1,\ldots,X_{n_i}]$ and isomorphisms $X_i\cong\spm(k<X_1,\ldots,X_{n_i}>/(f_{i1},\ldots,f_{ir_i}))$. To see this, let $Y=\spec(A)$ be an affine $k$-scheme of finite type, let $V\to Y^{\mbox{an}}$ be an open immersion of an affinoid rigid space $V$. Choose generators $f_1,\ldots, f_l$ of $A$ over $k$ such that $V\subset\{x\in Y^{\mbox{an}}|\quad f_i(x)\le 1 \mbox{ for all } i\}=Y_f$. By \cite{bgr},7.3.5/1, $V$ admits an admissible covering by rational subdomains of $Y_f$; by \cite{bgr},7.2.3/3, 7.2.4/1, the defining functions of these rational subdomains can be chosen in $A$.\\
Smooth rigid spaces are quasialgebraic (\cite{djvdp},3.4.1; \cite{elk},th.7).\\

\begin{satz} \label{raumkom} There exists a faithful functor $(.)'$ from the category of dagger spaces to the category of rigid spaces, together with a natural transformation $\nu:(.)'\to id(.)$ 
of functors with the following properties:\\(1) If $A'$ is the completion of the dagger algebra $A$, then $(\spm(A))'\cong\spm(A')$.\\
(2) $X$ is connected (resp. normal, resp. reduced, resp. regular) if and only if $X'$ is.\\(3) $\nu$ induces isomorphisms between the underlying $G$-topological spaces and between the local rings of the structure sheaves.\\
(4) $\rho:X\to Y$ is a closed immersion (resp. an open immersion, resp. a locally closed immersion, resp. an isomorphism, resp. quasicompact, resp. separated) if and only if $\rho':X'\to Y'$ is.\\(5) If $\rho$ is finite, so is $\rho'$.\\
(6) If $Y$ is a quasialgebraic rigid space, there exists an admissible affinoid covering $Y=\cup_iV_i$, affinoid dagger spaces $U_i$ and isomorphisms $U_i'\cong V_i$.\\ 
\end{satz}

{\sc Proof:} If $X$ is a dagger space, we have to construct a rigid structure sheaf on the underlying $G$-topological space of $X$. Locally this construction is prescribed by (1), and by the universal property of completion these local constructions glue. The other claims follow from our earlier observations.\\

\addtocounter{satz}{1}{\bf \arabic{section}.\arabic{satz}}  \newcounter{kohmodfu1}\newcounter{kohmodfu2}\setcounter{kohmodfu1}{\value{section}}\setcounter{kohmodfu2}{\value{satz}} If $X$ is a dagger space, then there is also a faithful functor $(.)'$ from the category of coherent 
${\cal O}_X$-modules to the category coherent ${\cal O}_{X'}$-modules, together with a natural transformation of functors $\mu:\nu^{-1}((.)_{Ab})\to((.)')_{Ab}$; 
here $\nu$ is from \ref{raumkom} and $(.)_{Ab}$ denotes the functor 'underlying abelian sheaf of a sheaf of modules'. 
Namely, if $X=\spm(A)$ is affinoid and $A'$ is the completion of $A$, and if ${\cal F}$ is a coherent ${\cal O}_X$-module, then ${\cal F}'={\cal O}_{X'}\otimes({\cal F}(X)\otimes_AA')$; this construction globalizes.\\
 
\addtocounter{satz}{1}{\bf \arabic{section}.\arabic{satz}}  For an open immersion $\spm(C)\subset\spm(D)$ of affinoid rigid spaces of the type $D\cong T_n(\rho)/I\to C\cong T_n/I.T_n$ for some $I$ and $\rho>1$ we write $\spm(C)\subset^{\dagger}\spm(D)$. 
It defines a structure of affinoid dagger space on the underlying $G$-topological space $|\spm(C)|$ of $\spm(C)$ (while in general an open immersion of the type $\spm(C)\subset\subset\spm(D)$ is only enough to establish a structure of dagger space on $|\spm(C)|$ 
which is not necessarily affinoid).\\
We also use the notation $\spm(A)\subset^{\dagger}\spm(B)$ for an open immersion of dagger spaces if $\spm(A')\subset^{\dagger}\spm(B')$ holds for the completions.\\

\addtocounter{satz}{1}{\bf \arabic{section}.\arabic{satz}}  Recall that if $U$ is an admissible open subset of the rigid space $W$, an admissible open subset $V\subset W$ is called a strict neighbourhood of $U$ in $W$ if $\{W-U,V\}$ is an admissible covering of $W$.\\
Given $\rho_0>1, I<T_n(\rho_0)$ and an open immersion of rigid spaces $\spm(T_n(\rho_0)/I)\to W$, the set $(V_{\rho}=\spm(T_n(\rho)/I.T_n(\rho)))_{\rho_0\ge\rho>1}$ is cofinal in the system of all strict neighbourhoods of $V_1=\spm(T_n/I.T_n)$ in $W$. To see this we may suppose $W=V_{\rho_0}$. 
If then $T$ is a strict neighbourhood of $V_1$ in $W$, the covering $W=\cup_{\rho>1}(W-V_{\rho})\cup T$ is admissible and therefore has a finite subcovering because $W$ is quasicompact. This implies $T\subset V_{\rho}$ for some $\rho>1$.\\

\addtocounter{satz}{1}{\bf \arabic{section}.\arabic{satz}}\newcounter{eindag1}\newcounter{eindag2}\setcounter{eindag1}{\value{section}}\setcounter{eindag2}{\value{satz}}  We deduce the following formula: Let $X$ be an affinoid dagger space, $X'\subset^{\dagger}Y$ an immersion of the type described above, ${\cal F}$ a coherent ${\cal O}_X$-module and 
${\cal F'}$ the associated coherent ${\cal O}_{X'}$-module. 
If ${\cal G}$ is a coherent ${\cal O}_{Y}$-module with ${\cal G}|_{X'}={\cal F}'$ (by \ref{mohomex} such a ${\cal G}$ always exists after perhaps shrinking $Y$), there is an identification $$\Gamma(X,{\cal F})=\{s\in\Gamma(X',{\cal G}) |\mbox{ there exists an}$$
$$ \mbox{extension of } s \mbox{ to a strict neighbourhood of } X' \mbox { in } Y\}.$$

\addtocounter{satz}{1}{\bf \arabic{section}.\arabic{satz}}  Let $X=\spm(A)\to Y=\spm(B)$ be a morphism of affinoid dagger spaces and $U\subset X$ an affinoid subdomain. We write $U\subset\subset_YX$ if there exists a surjection $\tau:B<X_1,\ldots,X_r>^{\dagger}\to A$ 
and $\epsilon\in\Gamma^*, \epsilon<1$ such that $U\subset \spm(A<\epsilon^{-1}.\tau(f_1),\ldots,\epsilon^{-1}.\tau(f_r)>^{\dagger})$.\\
A morphism $f:X\to Y$ of dagger (or rigid) spaces is called partially proper if $f$ is separated and there exist admissible affinoid coverings $Y=\cup Y_i$ and $f^{-1}(Y_i)=\cup X_{ij}$ (all $i$), 
such that for every $X_{ij}$ there is an admissible open affinoid subset $\tilde{X}_{ij}\subset f^{-1}(Y_i)$ with $X_{ij}\subset\subset_{Y_i}\tilde{X}_{ij}$ (cf.\cite{huet},p.59).\\
$f$ is called proper if it is quasicompact and partially proper.\\A dagger  (or rigid) space $X$ is called a Stein space if it admits an admissible affinoid covering $X=\cup_{i\in{\mathbb N}}U_i$ such that $U_i\subset^{\dagger}U_{i+1}$ for all $i$; we call $(U_i)_i$ a Stein covering.\\

\addtocounter{satz}{1}{\bf \arabic{section}.\arabic{satz}}  Stein spaces and spaces without boundary (\cite{lu},5.9) are partially proper. Compositions of partially proper morphisms of rigid spaces are partially proper (\cite{huet}). 
The (rigid) analytification of a $k$-scheme of finite type is partially proper. If $\rho:X\to Y$ is a partially proper morphism of dagger spaces, the associated morphism $\rho':X'\to Y'$ of rigid spaces is partially proper.\\
 
\begin{satz} \label{moduaeq} If $X$ is a partially proper dagger space, the functor $(.)':coh({\cal O}_X)\to coh({\cal O}_{X'})$ from \arabic{kohmodfu1}.\arabic{kohmodfu2} is an equivalence of categories. For a coherent ${\cal O}_X$-module ${\cal M}$ one has canonically $\Gamma(X,{\cal M})=\Gamma(X',{\cal M'})$.
\end{satz}

{\sc Proof:} (a) First consider the case where $X$ is a Stein space. Let $X=\cup_{i\in{\mathbb N}}U_i$ be a Stein covering and ${\cal G}$ a coherent ${\cal O}_{X'}$-module. 
The formula \arabic{eindag1}.\arabic{eindag2} defines for every $U_i$ a finite ${\cal O}_X(U_i)$-module ${\cal F}(U_i)$ and therefore a coherent ${\cal O}_X|_{U_i}$-module ${\cal F}|_{U_i}$. 
These glue to give a coherent ${\cal O}_X$-module ${\cal F}$ with ${\cal F}'={\cal G}$. Also ${\cal F}(X)={\cal G}(X)$ follows from \arabic{eindag1}.\arabic{eindag2}.\\
Now let ${\cal F}, {\cal G}$ be coherent ${\cal O}_X$-modules and $g:{\cal F}'\to{\cal G}'$ a morphism. $g$ induces morphisms ${\cal F}'(U_{i\rho})\to{\cal G}'(U_{i\rho})$ for all strict neighbourhoods $U_{i,\rho}$ 
of $U_i$ in $U_{i+1}$, and by \arabic{eindag1}.\arabic{eindag2} this defines ${\cal O}_X(U_i)$-linear morphisms ${\cal F}(U_i)\to{\cal G}(U_i)$, that is ${\cal O}_X|_{U_i}$-linear morphisms ${\cal F}|_{U_i}\to{\cal G}|_{U_i}$. These glue to give $f:{\cal F}\to{\cal G}$ with $f'=g$.\\
(b) For the general case it is now enough to remark that $X$ has an admissible covering $X=\cup_jS_j$ such that all $S_j$ and their finite intersections are Stein spaces. For example, start with an affinoid covering $X=\cup_jV_j$ such that for all $V_j$ there is an admissible open affinoid subspace $Z_j\subset X$ with $V_j\subset^{\dagger}Z_j$. 
Take $S_j$ to be the interior of $Z_j$ (here by the interior of an affinoid dagger or rigid space $T=\spm(A)$ with respect to a fixed closed immersion $T\hookrightarrow{\bf D}^m$ we mean $T\cap({\bf D}^0)^m$; for radii of polydiscs different from $1$ one has obvious variants). By separatedness of $X$ all finite intersections of the $S_j$ again admit a closed immersion into a polydisc (of some polyradius) without boundary, in particular are Stein spaces.\\

\begin{satz}\label{raumaeq} The functor $(.)'$ from \ref{raumkom} induces an equivalence between the category of partially proper dagger spaces and the category of partially proper rigid spaces.\\ 
\end{satz}

{\sc Proof:} $(.)'$ is essentially surjective: Here for any partially proper rigid space $Z$ we have to give a sheaf ${\cal R}$ of $k$-algebras which at open  affinoids $U_1\subset Z$, 
for which $U_2\subset Z$ exist with $U_1\subset^{\dagger}U_2$, has as value a dagger algebra ${\cal R}(U_1)$ whose completion is ${\cal O}_Z(U_1)$. This is done by means of the construction in \ref{moduaeq}.\\
$(.)'$ is full: Let $X, Y$ be partially proper dagger spaces and $g:X'\to Y'$ a morphism. Choose a family $\{U_i^{\delta_i}\}_i$ of open subspaces of $Y'$ of the following type: 
$U_i^{\delta_i}=\spm(T_{n_i}(\delta_i)/I_i)$ for some $\delta_i>1, n_i\in{\mathbb N}$ and $I_i<T_{n_i}(\delta_i)$, and if $U_i^{\delta}=\spm(T_{n_i}(\delta)/I_i.T_{n_i}(\delta))$ for $\delta_i\ge\delta\ge1$ the family $(U_i^1)_i$ 
forms an admissible covering of $Y'$. Furthermore write $U_i^{\delta,0}=\cup_{\tilde{\delta}<\delta}\spm(T_{n_i}(\tilde{\delta})/I_i.T_{n_i}(\tilde{\delta}))$, the interior of $U_i^{\delta}$. 
Identifying the underlying $G$-topological spaces the sets $g^{-1}(U_i^1)$ define open dagger subspaces $V_i$ of $X$, and the sets $U_i^1$ define open affinoid dagger subspaces $W_i$ of $Y$. 
Now it is enough to give compatible morphisms $f_i:V_i\to W_i$ of dagger spaces such that $f_i'=g|_{g^{-1}(U_i^1)}$. 
Define $f_i$ by giving a $k$-algebra morphism ${\cal O}_Y(W_i)\to{\cal O}_X(V_i)$ as follows: $g$ defines for all $\delta_i\ge\delta>1$ morphisms 
${\cal O}_{Y'}(U_i^{\delta})\to{\cal O}_{X'}(g^{-1}(U_i^{\delta}))$, and in the limit
$${\cal O}_Y(W_i)=\lim_{\stackrel{\to}{\delta\to 1}}{\cal O}_{Y'}(U_i^{\delta})\to \lim_{\stackrel{\to}{\delta\to 1}}{\cal O}_{X'}(g^{-1}(U_i^{\delta})).$$ 
On the other hand there is a canonical map $$\lim_{\stackrel{\to}{\delta\to 1}}{\cal O}_{X'}(g^{-1}(U_i^{\delta}))\to{\cal O}_X(V_i),$$
namely: $g^{-1}(U_i^{\delta,0})$ is partially proper for all $\delta_i\ge\delta>1$ (because $g$ is partially proper by \cite{bgr},9.6.2/4, and compositions of partially proper morphisms are partially proper). 
Therefore ${\cal O}_X(\rho'^{-1}(U_i^{\delta,0}))\to{\cal O}_{X'}(g^{-1}(U_i^{\delta,0}))$ is an isomorphism (cf. \ref{moduaeq}), and one gets $${\cal O}_{X'}(g^{-1}(U_i^{\delta}))\to{\cal O}_{X'}(g^{-1}(U_i^{\delta,0}))\cong{\cal O}_X(g^{-1}(U_i^{\delta,0}))\to{\cal O}_X(V_i).$$

\section{Cohomology}

\begin{pro}\label{affkoh} If $X=\spm(A)$ is an affinoid dagger space, every coherent ${\cal O}_X$-module $\cal{F}$ is generated by its global sections, and one has $H^n(X,{\cal F})=0$ for all $n\ge1$.
\end{pro}  

{\sc Proof:} Follows from \ref{kiehl} and \ref{modazyk}.\\

\begin{satz}\label{modukohaeq}\label{steiazyk} For a partially proper dagger space $X$, a coherent ${\cal O}_{X}$-module $\cal{F}$ and $i\ge0$ one has $H^i(X,{\cal F})=H^i(X',{\cal F'})$. If $X$ is a Stein space, these groups vanish for $i>0$.
\end{satz}

{\sc Proof:} First consider the case where $X$ is a Stein space. Let $X=\cup_{i\in{\mathbb N}}U_i$ be a Stein covering. 
As it is affinoid it is acyclic for ${\cal F}$ and induces an affinoid, hence ${\cal F}'$-acyclic covering $X'=\cup_{i\in{\mathbb N}}U'_i$. 
The assertion in case $i=0$ follows from \arabic{eindag1}.\arabic{eindag2}, the assertion in case $i>1$ is a formal consequence of the existence of nested acyclic coverings indexed by ${\mathbb N}$.\\
From \cite{kiaub} we infer $H^1(X',{\cal F}')=0$. Now let $(g_i)_{i\in{\mathbb N}}\in{\cal F}(U_i)_{i\in{\mathbb N}}$ be a cocycle. 
For every $i\in{\mathbb N}$ there exists an open affinoid $V_i\subset X'$ with $U_i'\subset^{\dagger}V_i\subset U_{i+1}'$ and a $g_i'\in{\cal F}'(V_i)$ inducing $g_i$ (see \arabic{eindag1}.\arabic{eindag2}). 
Because of $H^1(X',{\cal F}')=0$ the cocycle $(g_i')_{i\in{\mathbb N}}\in{\cal F}'(V_i)_{i\in{\mathbb N}}$ for the ${\cal F'}$-acyclic covering $X'=\cup_{i\in{\mathbb N}}V_i$ is a coboundary, 
that is there are $f_i'\in{\cal F'}(V_i)$ satisfying $f_{i+1}'-f_i'=g_i'$. 
The $f_i'\in{\cal F'}(V_i)$ induce elements $f_i\in{\cal F}(U_i)$ (again \arabic{eindag1}.\arabic{eindag2}) which testify that $(g_i)_{i\in{\mathbb N}}$ is a coboundary.\\
The case of a general $X$ is now deduced from the above by means of the \v{C}ech-spectral sequence associated with an admissible covering $X=\cup_jW_j$ such that all $W_j$ and their finite intersections are Stein spaces (as in the proof of \ref{moduaeq}).\\

\addtocounter{satz}{1}{\bf \arabic{section}.\arabic{satz}}  Besides the rigid analytification functor there is also a dagger analytification functor for $k$-schemes of finite type (to be defined in the same manner, or equivalently via \ref{raumaeq}). 
As a corollary of \ref{modukohaeq} we get that for proper $k$-schemes of finite type dagger analytification satisfies the GAGA-principle (since this is true for rigid analytification, \cite{kiend}).\\
Another corollary of \ref{modukohaeq} will be that for a smooth partially proper dagger space $X$ the de Rham cohomology of $X$ coincides with the de Rham cohomology of $X'$ 
(compare the spectral sequences $E^{pq}_1=H^q(\Omega^p)\Rightarrow H^{p+q}(\Omega^{\bullet})$ for $X$ and $X'$).\\ 

\addtocounter{satz}{1}{\bf \arabic{section}.\arabic{satz}}  We recall some notations from \cite{berco}. Suppose $k=\q(R)$ for a complete discrete valuation ring $R$ (of mixed characteristic). If ${\cal X}$ is an admissible formal $R$-scheme (\cite{bolu}) with generic fibre $X$ and specialization morphism 
$s:X\to{\cal X}$, and if $V\subset{\cal X}_s$ is a subset of its special fibre, we write $]V[=s^{-1}(V)$. 
If $j:V\to{\cal X}_s$ is an open immersion, the functor $j^{\dagger}$ from the category of abelian sheaves on $X$ to itself is defined by ${\cal F}\mapsto j^{\dagger}{\cal F}=\lim_{\stackrel{\to}{U}}j_{U*}j_U^{-1}{\cal F}$, 
where $j_U:U\to X$ runs through the strict neighbourhoods of $]V[$ in $X$.\\

\begin{satz}\label{kohaer} Let $X\stackrel{f}{\to}Y$ be a proper morphism of dagger spaces and ${\cal F}$ a coherent ${\cal O}_X$-module. 
Then $R^qf_*{\cal F}$ is a coherent ${\cal O}_Y$-module for all $q\in{\mathbb N}$ satisfying $(R^qf_*{\cal F})'\cong R^qf_*'({\cal F}')$, at least under the following two assumptions:\\
(1) $k$ is the fractionfield of a complete discrete valuation ring $R$ (of mixed characteristic), and\\
(2) there exists an admissible affinoid covering $Y=\cup_{i\in I}U_i$ with the following property: 
Setting $X\times_YU_i\stackrel{f_i}{\to}U_i$ there exists for all $i$ a proper morphism $\tilde{X}_i\stackrel{\tilde{f}_i}{\to}\tilde{U}_i$ of rigid spaces and a coherent ${\cal O}_{\tilde{X}_i}$-module $\tilde{\cal F}_i$ 
such that $(f_i',{\cal F}'|_{X_i'})$ arises from $(\tilde{f}_i,\tilde{\cal F}_i)$ by a base change of type $U_i'\subset^{\dagger}\tilde{U}_i$.
\end{satz}

{\sc Proof:} We may suppose $Y=\spm(A)$ and $(f',{\cal F}')$ arises from $(\tilde{X}\stackrel{\tilde{f}}{\to}\tilde{Y},\tilde{\cal F})$ by the base change $Y'\subset^{\dagger}\tilde{Y}=\spm(\tilde{A})$ (where $\tilde{f}$ is proper and $\tilde{\cal F}$ a coherent ${\cal O}_{\tilde{X}}$-module). 
Choose a finite affinoid covering $X=\cup_{i\in K}V_i$ such that for every $i\in K$ there is an open affinoid $\tilde{W}_i\subset\tilde{X}$ with $V_i'\subset^{\dagger}\tilde{W}_i$ (think of $X'$ 
as an open subspace of $\tilde{X}$). 
Then choose an admissible formal $R$-scheme $\tilde{\cal X}$ with generic fibre $\tilde{X}$ such that for every $i\in K$ there is an open immersion $Z_i\to\tilde{\cal X}_s$ into the special fibre with $]Z_i[_{\tilde{\cal X}}=V_i'$ (\cite{bolu}). 
Set $Z=\cup_{i\in K}Z_i\stackrel{j}{\to}\tilde{X}_s$, that is $]Z[_{\tilde{\cal X}}=X'$. 
The open immersion $X'\stackrel{l}{\to}\tilde{X}$ is acyclic for ${\cal F}$, hence $H^q(X,{\cal F})=H^q(\tilde{X},l_*{\cal F})$ for all $q$. 
There is an obvious morphism of sheaves $l_*{\cal F}\to j^{\dagger}\tilde{\cal F}$ inducing maps $H^q(X,{\cal F})=H^q(\tilde{X},l_*{\cal F})\to H^q(\tilde{X},j^{\dagger}\tilde{\cal F})$. Claim: 
These are isomorphisms. By \cite{berco},2.1, $H^q(\tilde{X},j^{\dagger}\tilde{\cal F})$ can be computed locally with respect to the covering $Z=\cup_{i\in K}Z_i$, 
we therefore have to show $H^q(V_J,{\cal F})\cong H^q(\tilde{X},j_J^{\dagger}\tilde{\cal F})$ for all $q$, all $J\subset K$, where we put $V_J=\cap_{i\in J}V_i$ and $(\cap_{i\in J}Z_i)=Z_J\stackrel{j_J}{\to}\tilde{\cal X}_s$. 
For $q=0$ this is clear, for $q>0$ one has $H^q(V_J,{\cal F})=0$ (because $V_J$ is affinoid), but also $H^q(\tilde{X},j_J^{\dagger}\tilde{\cal F})=0$ (formation of cohomology commutes with direct limits because $\tilde{X}$ 
is quasicompact; $V_J'=]Z_J[_{\tilde{\cal X}}$ has a cofinal system of strict neighbourhoods in $\tilde{X}$ consisting of affinoids). 
The claim follows.\\Now let $(Y_{\rho})_{\rho>1}$ be a system of neighbourhoods of $Y'$ in $\tilde{Y}$ of type $Y'\subset^{\dagger}Y_{\rho}=\spm(A_{\rho})$ such that $\lim_{\stackrel{\leftarrow}{\rho\to 1}}Y_{\rho}=Y$. 
Set $Y_{\rho}\times_{\tilde{Y}}\tilde{X}=X_{\rho}\stackrel{l_{\rho}}{\to}\tilde{X}$. Then
$$H^q(\tilde{X},j^{\dagger}\tilde{\cal F})\cong\lim_{\stackrel{\to}{\rho\to 1}}H^q(\tilde{X},l_{\rho *}(\tilde{\cal F}|_{X_{\rho}}))\cong\lim_{\stackrel{\to}{\rho\to 1}}H^q(X_{\rho},\tilde{\cal F})$$
$$\stackrel{(*)}{\cong}\lim_{\stackrel{\to}{\rho\to 1}}(H^q(\tilde{X},\tilde{\cal F})\otimes_{\tilde{A}}A_{\rho})=H^q(\tilde{X},\tilde{\cal F})\otimes_{\tilde{A}}A.$$
Here we apply \cite{kiend} to get $(*)$. Together we obtain $H^q(X,{\cal F})=H^q(\tilde{X},\tilde{\cal F})\otimes_{\tilde{A}}A$ for the finite $\tilde{A}$-module $H^q(\tilde{X},\tilde{\cal F})$. 
If $U\subset Y$ is open affinoid such that there is an open $\tilde{U}\subset\tilde{Y}$ with $U'\subset^{\dagger}\tilde{U}$ the same can be shown for $U$ instead of $Y$. The theorem follows.\\

\addtocounter{satz}{1}{\bf \arabic{section}.\arabic{satz}}  Certainly both assumptions (1) and (2) are superfluous, perhaps (2) is even automatic.\\

\section{Duality and K\"unneth formula}
\addtocounter{satz}{1}{\bf \arabic{section}.\arabic{satz}}  For a dagger algebra $A$ we define $d:A\to\Omega_A$ to be the universal $k$-linear derivation of $A$ into finite $A$-modules. 
In the usual way we get for a smooth dagger space $X$ a $k$-linear complex $\Omega_X^{\bullet}$ with all $\Omega_X^i$ being coherent locally free ${\cal O}_X$-modules such that $(\Omega_X^i)'=\Omega_{X'}^i$ 
(smoothness of $X$ is defined as for rigid spaces, or equivalently by requiring that $X'$ be smooth). In particular, put $\omega_X=\Omega_X^n$ if $X$ is of pure dimension $n$. \\

\addtocounter{satz}{1}{\bf \arabic{section}.\arabic{satz}} We assume from now on that $k$ is spherically complete. We endow dagger algebras with a topology which is finer than the normtopology from \ref{normvgl} as follows: 
If $A\cong\lim_{\stackrel{\to}{\rho\to 1}}T_n(\rho)/I.T_n(\rho)$ with $\rho_0\ge\rho>1$ 
for some $\rho_0>1$ and $I<T_n(\rho_0)$, we define the direct limit topology on $A$ to be the finest locally $k$-convex topology such that all maps $T_n(\rho)/I.T_n(\rho)\to A$ are continuous. 
The direct limit topology on a finite $A$-module $M$ is the quotienttopology of the direct limit topology of $A^r$ with respect to a surjection $A^r\to M$ (some $r\in{\mathbb N}$). 
Using \ref{extens} we see that these definitions are independent of the chosen representations, and from \cite{mori},3.3, 3.5, we deduce that in this way $A$ becomes a complete reflexive Hausdorffspace. In the following we will only consider these topologies.\\

\addtocounter{satz}{1}{\bf \arabic{section}.\arabic{satz}}  For an affinoid dagger space $X=\spm(A)$ and a coherent ${\cal O}_X$-module ${\cal F}$ we define $H_c^*(X,{\cal F})$ as follows. Choose a representation $A=W_m/I.W_m$ for some $\rho_0>1, I<T_m(\rho_0)$ 
and set $A_{\rho}=T_m(\rho)/I.T_m(\rho)$ and $X_{\rho}=\spm(A_{\rho})$ for $1<\rho\le\rho_0$. 
After shrinking $\rho_0$ if necessary, there is a coherent ${\cal O}_{X_{\rho_0}}$-module ${\cal F}_{\rho_0}$ such that ${\cal F}(X)=\lim_{\stackrel{\to}{\rho>1}}{\cal F}_{\rho_0}(X_{\rho})$. 
Put $H_c^i(X,{\cal F})=H_X^i({\cal F}_{\rho_0})$ for $i\ge0$. 
Endow $H_c^i(X,{\cal F})=H_X^i({\cal F}_{\rho})$ for $i>0$ with the finest locally $k$-convex topology such that all boundary maps  $H^{i-1}(X_{\rho_0}-X,{\cal F}_{\rho_0})\to H_X^i({\cal F}_{\rho_0})$ are continuos (here $H^{i-1}(X_{\rho}-X,{\cal F}_{\rho})$ carries the topology described in \cite{vdpse},1.6,\cite{bey}). 
Note that $H_c^i(X,{\cal F})$ carries also a natural structure of topological $A$-module.\\

\begin{satz}\label{serduarec}\label{sedu} Let $X$ be a smooth affinoid dagger space of pure dimension $d$ and ${\cal F}$ a coherent ${\cal O}_X$-module. 
Then for all $0\le i\le d$ there is a pairing$$\ext_{{\cal O}_X}^i({\cal F},\omega_X)\times H_c^{d-i}(X,{\cal F})\to k,$$ functorial in ${\cal F}$ and inducing isomorphisms 
$$\ext_{{\cal O}_X}^i({\cal F},\omega_X)\cong H_c^{d-i}(X,{\cal F}){\bf {\bf \check{}}}=\ho_{k,cont}(H_c^{d-i}(X,{\cal F}),k)$$ and 
$$\quad H_c^{d-i}(X,{\cal F})\cong\ext_{{\cal O}_X}^i({\cal F},\omega_X){\bf \check{}}=\ho_{k,cont}(\ext_{{\cal O}_X}^i({\cal F},\omega_X),k).$$
\end{satz}

{\sc Proof:} We begin by considering the case $X=\spm(A)$ with $A=k<T_1,\ldots,T_m>^{\dagger}$. For $\rho>1$ put $A_{\rho}=k<\rho^{-1}.T_1,\ldots,\rho^{-1}.T_m>$ and $X_{\rho}=\spm(A_{\rho})$. 
Then there is a canonical identification (compare \cite{etst},\cite{vdpse},\cite{bey})$$\quad H_c^m(X,\omega_X)=\{\sum_{\stackrel{\mu\in{\mathbb Z}^m}{\mu<0}}a_{\mu}T^{\mu}dT |\quad \eta^{-|\mu|}|a_{\mu}|\stackrel{|\mu|\to\infty}{\to}0\mbox{ for all }\eta>1\}.$$
Arguing as in \cite{rocri},6.5 (where the case $m=1$ is done; compare also \cite{etst}) we see that the pairing$$A\times H_c^m(X,\omega_X)\stackrel{\res}{\to}k,$$
$$ (\sum_{\stackrel{\alpha\in{\mathbb Z}^m}{\alpha\ge0}}b_{\alpha}T^{\alpha},\sum_{\stackrel{\mu\in{\mathbb Z}^m}{\mu<0}}a_{\mu}T^{\mu}dT)\mapsto \sum_{\stackrel{\alpha\in{\mathbb Z}^m}{\alpha\ge0}}b_{\alpha}a_{-\alpha-1}$$
induces isomorphisms $A\cong\ho_{k,cont}(H_c^m(X,\omega_X),k)$ and $H_c^m(X,\omega_X)\cong\ho_{k,cont}(A,k)$. From this we deduce the general case by standard arguments (similar to those in \cite{vdpse},\cite{bey}), analysing a closed immersion $X\to\spm(W_m)$.\\

\addtocounter{satz}{1}{\bf \arabic{section}.\arabic{satz}} \newcounter{mehrsedu1}\newcounter{mehrsedu2}\setcounter{mehrsedu1}{\value{section}}\setcounter{mehrsedu2}{\value{satz}} From \ref{sedu} one gets Serre duality for more general quasicompact dagger spaces. Also, passing to the limit in a Stein covering, one gets the wellknown Serre duality on Stein spaces.\\

\begin{lem}\label{dualstet} Let $K^{\bullet}=(K^i\stackrel{d^i}{\to}K^{i+1})_{i\in{\mathbb Z}}$ be a complex of locally $k$-convex complete vector spaces, all $d^i$ being strict. 
Setting $L^{\bullet}=\ho_{k,cont}(K^{\bullet},k)$ one has $H^q(L^{\bullet})=\ho_{k,cont}(H^q(K^{\bullet}),k)$ for all $q$.\\
\end{lem}

{\sc Proof:} Recall that a morphism $\phi:E\to F$ of topological groups is called strict, if $\phi:E\to\phi(E)$ is open, where $\phi(E)$ carries the topology induced from $F$. We observe that \ref{dualstet} is an easy consequence of the Hahn-Banach theorem: 
If $E$ is locally $k$-convex complete $k$-vector space, $M\subset E$ a subspace and $M\stackrel{\lambda}{\to}k$ a continuous linear form, then $\lambda$ extends to a continuous linear form on $E$. 
For the proof of this we refer to \cite{roo},thm.4.10, in case $E$ is a Banach space, to which the general case is reduced as usual.\\ 

\begin{lem}\label{diffstri} Suppose that $k$ is discretely valued and of characteristic $0$.\footnote{Unfortunately, in the published version of the present paper these hypotheses are missing. I thank Mark Kisin for pointing out that they should be imposed here in order to have the finiteness result of \cite{en2dag} available.} Let $X$ be a smooth Stein space or a smooth affinoid dagger space. All differentials $\Omega_X^i(X)\stackrel{d_X^i}{\to}\Omega_X^{i+1}(X)$ are strict and have closed image.
\end{lem}

{\sc Proof:} The second claim follows from the first by the completeness of $\Omega_X^i(X)$. For the first we begin with the case where $X$ is Stein. By \cite{lustein} there is a closed immersion $X\to{\bf A}^n$ into the (analytification of the) affine space. 
The canonical surjections $\Omega^j_{{\bf A}^n}({\bf A}^n)\stackrel{p_j}{\to}\Omega^j_X(X)$ are strict, and since we have $p_{i+1}\circ d^i_{\bf{A}^n}=d^i_X\circ p_i$ 
we are therefore reduced to the case $X={\bf A}^n$. 
In this case $\Omega_X^{\bullet}(X)$ is acyclic in positive degrees, hence $\bi(d_X^i)=\ke(d_X^{i+1})$ is closed in $\Omega_X^{i+1}(X)$. But this means that $\bi(d_X^i)$ is a Fr\'{e}chet space 
(for its induced topology from $\Omega_X^{i+1}(X)$). By the theorem of Banach, $\Omega_X^i(X)\stackrel{d_X^i}{\to}\bi(d_X^i)$ is open.\\ 
Now we consider the case where $X=\spm(A)$ is an affinoid dagger space. As above we reduce to the case $A=W_n$. For $X_{\rho}=\lim_{\stackrel{\to}{\rho'<\rho}}\spm(T_n(\rho'))$ 
we endow $\Omega^j_{X_{\rho}}(X_{\rho})=\lim_{\stackrel{\leftarrow}{\rho'<\rho}}\Omega^j_{T_n(\rho')}$ 
with its inverse limit topology and note that the direct limit topology on $\Omega_A^j$ 
as defined above is also the finest locally $k$-convex topology on $\Omega_A^j=\lim_{\to}\Omega^j_{X_{\rho}}(X_{\rho})$ such that all canonical maps $\Omega^j_{X_{\rho}}(X_{\rho})\to\Omega^j_A$ are continuous. 
Therefore it suffices to show strictness of all maps $\Omega^j_{X_{\rho}}(X_{\rho})\stackrel{d^j_{\rho}}{\to}\Omega^{j+1}_{X_{\rho}}(X_{\rho})$, which just has been done, because the $X_{\rho}$ are Stein spaces.\\

\addtocounter{satz}{1}{\bf \arabic{section}.\arabic{satz}} From now on we assume that $k$ is discretely valued and $\kara(k)=0$. For a smooth dagger space $X$ we define $H_{dR}^*(X)=H^*(X,\Omega_X^{\bullet})$. If $X=\spm(W_m/I.W_m)$ for some $\rho>1, I<T_m(\rho)$ we define  $H_{dR,c}^*(X)=H^*_X(\Omega^{\bullet}_{\spm(T_m(\rho)/I.T_m(\rho))})$ (cohomology with support in the closed subset $X\subset\spm(T_m(\rho)/I.T_m(\rho))$).\\

\begin{satz}\label{dagpoin} If $X$ is a smooth affinoid dagger space of pure dimension $n$, there are canonical isomorphisms $H_{dR}^i(X)\cong \ho_{k,cont}(H_{dR,c}^{2n-i}(X),k)$ and $H_{dR,c}^{2n-i}(X)=\ho_{k,cont}(H_{dR}^i(X),k)$ for all $0\le i\le 2n$.
\end{satz}

{\sc Proof:} We begin by observing that for all locally free coherent ${\cal O}_X$-modules ${\cal F}$ we have $H_c^q(X,{\cal F})=0$ for all $q\ne n$. Arguing along a finite surjective map $f:X\to\spm(W_n)$ 
and an extension of $f$ to a finite surjective map over some $\spm(T_n(\rho))$, this statement can be reduced to the case $X=\spm(W_n)$ and ${\cal F}={\cal O}_X$ and can then be checked by reasoning as in \cite{vdpse}, \cite{bey}.\\ 
This implies $H_{dR,c}^{q+n}(X)=H^q(H^n_c(\Omega^{\bullet}_X))$ for all $q\ge0$, where we denote by $H_c^n(\Omega^{\bullet}_X)$ the complex$$\ldots\to H_c^n(X,\Omega^i_X)\to H_c^n(X,\Omega_X^{i+1})\to\ldots.$$
But this complex is canonically identified with $$\ldots\to H_c^n(X,{\cal O}_X)\otimes_A\Omega^i_A\to H_c^n(X,{\cal O}_X)\otimes_A\Omega^{i+1}_A\to\ldots$$
and now the theorem follows by combining \ref{sedu}, \ref{dualstet} and \ref{diffstri}.\\

\addtocounter{satz}{1}{\bf \arabic{section}.\arabic{satz}}  The $k$-vector spaces $H_{dR}^i(X)$ are in fact finite dimensional if $k$ is discretely valued (\cite{doktor},\cite{en2dag}). By \ref{diffstri} they are topologically separated, therefore in this case \ref{dagpoin} becomes an algebraic duality between finite dimensional vector spaces.\\ 

\begin{satz} If $X$ is a smooth Stein space of pure dimension $n$, there are canonical isomorphisms $H_{dR}^i(X)\cong \ho_{k,cont}(H_{dR,c}^{2n-i}(X),k)$ and $H_{dR,c}^{2n-i}(X)=\ho_{k,cont}(H_{dR}^i(X),k)$ for all $0\le i\le 2n$.
\end{satz}

{\sc Proof:} Let $X=\cup_{j\in{\mathbb N}}U_j$ be an admissible covering by affinoid dagger spaces. The proof of \ref{dagpoin} gives $H_c^i(X,{\cal F})\cong \lim_{\stackrel{\to}{j}}H_{U_j}^i(U_{j+1},{\cal F}|_{U_{j+1}})=0$ for all $i\ne n$, all coherent ${\cal O}_X$-modules ${\cal F}$. 
So we can conclude as in \ref{dagpoin}, this time combining Serre duality of Stein spaces (cf. \arabic{mehrsedu1}.\arabic{mehrsedu2}, \cite{vdpse},\cite{bey}) with \ref{dualstet} and \ref{diffstri}.\\

\begin{satz}\label{kunneth} Let $k$ be discretely valued and let $X$ and $Y$ be smooth dagger spaces. There are canonical isomorphisms $\oplus_{n=p+q}H_{dR}^p(X)\otimes_kH_{dR}^q(Y)\cong H_{dR}^n(X\times Y)$.
\end{satz}

{\sc Proof:} Since $X\times Y$ can be admissibly covered by subspaces of type $U\times V$, where $U\subset X$ and $V\subset Y$ are open affinoid, we may assume $X$ and $Y$ affinoid, $X=\spm(A)$ and $Y=\spm(B)$. We have to show that if $M^{\bullet}=(M^{\bullet},d^{\bullet}_M)$ is the associated simple complex of the double-complex $\Omega_A^{\bullet}\otimes_k\Omega_B^{\bullet}$, the canonical morphism $(M^{\bullet},d^{\bullet}_M)\to(\Omega_{A\otimes_k^{\dagger}B}^{\bullet},d^{\bullet}_{\otimes^{\dagger}})$ is a quasiisomorphism. Note that this is a morphism of complexes of topological $k$-vector spaces. One verifies that all maps $M^j\to\Omega_{A\otimes_k^{\dagger}B}^j$ are injective, strict and have a dense image. By \ref{diffstri} the cohomology vector spaces $H^n(M^{\bullet})$ are separated, i.e. $\bi(d_M^{n-1})$ is closed in $M^n$, implying that the maps $H^n(M^{\bullet})\to H^n(\Omega_{A\otimes_k^{\dagger}B}^{\bullet})$ are injective and strict. But $H^n(M^{\bullet})$ is finite dimensional by \cite{doktor},\cite{en2dag}, hence complete. Therefore all we must show is that $\ke(d_M^n)$ is dense in $\ke(d^n_{\otimes^{\dagger}})$.\\
Write $A=W_{n_1}/(f_1)$, $B=W_{n_2}/(f_2)$ with tuples $f_i\in T_{n_i}(\rho_0)^{r_i}$ for some $\rho_0>1, r_i\in{\mathbb N}$ (here $i=1,2$). For $1<\rho\le\rho_0$ set $A_{\rho}=T_{n_1}(\rho)/(f_1),\quad\quad B_{\rho}=T_{n_2}(\rho)/(f_2),\quad\quad X_{\rho}=\spm(A_{\rho}),\quad\quad Y_{\rho}=\spm(B_{\rho})$ and $Z_{\rho}=\spm(A_{\rho}\hat{\otimes}_kB_{\rho})$. Let $X^0_{\rho}=\cup_{1<\rho'<\rho}X_{\rho'},\quad\quad Y^0_{\rho}=\cup_{1<\rho'<\rho}Y_{\rho'}$ and $Z^0_{\rho}=\cup_{1<\rho'<\rho}Z_{\rho'}$ be the interiors. Consider the differentials$$\oplus_{n=p+q}\Omega^p_{X^0_{\rho}}(X^0_{\rho})\otimes_k\Omega^q_{Y^0_{\rho}}(Y^0_{\rho})\stackrel{d^n_{\rho}}{\to}\oplus_{n+1=p+q}\Omega^p_{X^0_{\rho}}(X^0_{\rho})\otimes_k\Omega^q_{Y^0_{\rho}}(Y^0_{\rho}).$$Using strictness (\ref{diffstri}) and the fact that the topologies on the spaces $\Omega^n_{Z^0_{\rho}}(Z^0_{\rho})$ have a countable basis (they are Fr\'{e}chet spaces), we get (as in the proof of \cite{bgr},1.1.9/5) that $\ke(d^n_{\rho})$ is dense in $\ke(\Omega^n_{\tilde{Z}_{\rho}}(\tilde{Z}_{\rho})\stackrel{\hat{d}^n_{\rho}}{\to}\Omega^{n+1}_{\tilde{Z}_{\rho}}(\tilde{Z}_{\rho}))$. Passing to the limit as $\rho\to 1$, we conclude.\\

\section{Comparison with Berthelot's rigid cohomology}

We assume $k=\q(R)$ for a complete discrete valuation ring $R$ of mixed characteristic.\\

\begin{satz} \label{elmber} Let ${\cal X}$ be a proper admissible (\cite{bolu}) formal $\spf(R)$-scheme, $Y\to{\cal X}_s$ an immersion into its special fibre with schematic closure $j:Y\to\bar{Y}$ in ${\cal X}_s$. 
Then $]\bar{Y}[_{\cal X}$ is a partially proper rigid space, therefore equivalent with a dagger space $Q$. Let $X$ be the open subspace of $Q$ whose underlying set is identified with $]Y[_{\cal X}$. Let $q\in{\mathbb N}$.\\
(a) If ${\cal F}$ is a coherent ${\cal O}_Q$-module, ${\cal F}'$ the associated coherent ${\cal O}_{]\bar{Y}[_{\cal X}}$-module, then there is a canonical isomorphism $H^q(X,{\cal F}_X)\cong H^q(]\bar{Y}[_{\cal X},j^{\dagger}{\cal F}').$\\
(b) If $X$ is smooth, there is a canonical isomorphism $H^q_{dR}(X)\cong H^q(]\bar{Y}[_{\cal X},j^{\dagger}\Omega^{\bullet}_{]\bar{Y}[_{\cal X}})$.\\
(c) If ${\cal X}$ is smooth along $Y$, there is a canonical isomorphism
$H^q_{dR}(X)\cong H_{rig}^q(Y/k)$.
\end{satz}

{\sc Proof:} (b) follows from (a), and (c) follows from (b) by the definition of rigid cohomology (\cite{berco}).\\
Let $Y=\cup_{i\in I}Y_i$ be an affine open covering, for $\emptyset\ne J\subset I$ let $Y_J=\cap_{i\in J}Y_i\stackrel{j_J}{\to}\bar{Y}$ and $]Y_J[_{\cal X}\stackrel{i_J}{\to}]\bar{Y}[_{\cal X}$ be the open immersions. 
Since all $i_J$ are affinoid, the Cechcomplex $L^{\bullet}$ built out of all $i_{J*}{\cal F}_{]Y_J[_{\cal X}}$ resolves $Ri_*{\cal F}_X$ (where $]Y[_{\cal X}\stackrel{i}{\to}]\bar{Y}[_{\cal X}$ is the open immersion, and the underlying $G$-topological spaces of $Q$ and $]\bar{Y}[_{\cal X}$ are identified). On the other hand, \cite{berco},2.1, tells us that the Cechcomplex $K^{\bullet}$ built out of all $j_J^{\dagger}{\cal F}'$ resolves $j^{\dagger}{\cal F}'$. There is a canonical map $L^{\bullet}\to K^{\bullet}$. 
Thus to prove (a) we may assume $Y$ is affine. We may then also assume that there is an open affinoid $U\subset{\cal X}_k$ and an isomorphism $U=\spm(T_n(\rho')/I)$ for some $1<\rho'\in\Gamma^*, I<T_n(\rho')$ 
such that $]W[_{\cal X}=\spm(T_n/I.T_n)\subset U$ for an open affine $W\subset{\cal X}_s$ for which $Y\to{\cal X}_s$ factorizes via a closed immersion $Y\to W$ (such a $U$ exists, at least after further refinement of the $Y$-covering, use \cite{lu}). 
We may finally assume that there are $f_1,\ldots,f_r\in T_n(\rho')/I$ such that $U_{\bar{Y}}=]\bar{Y}[_{\cal X}\cap U=\{x\in U|\quad |f_i(x)|<1\mbox{ for all }i=1,\ldots,r\}$. 
Then for $\rho,\mu\in\Gamma^*$ with $\mu<1$ and $1<\rho\le\rho'$ set$$U_{\mu,\rho}=\spm((T_n(\rho)/I.T_n(\rho))<\mu^{-1}.f_1,\ldots,\mu^{-1}.f_r>),$$
$$V_{\mu}=\spm((W_n/I.W_n)<\mu^{-1}.f_1,\ldots,\mu^{-1}.f_r>^{\dagger})$$ (viewed as subspaces of ${\cal X}_k$ resp. $Q$). Let $(\mu_n)_{n\in{\mathbb N}}\to 1$ be a monotonicly increasing sequence in $\Gamma^*$. 
The covering $X=\cup_{n\in{\mathbb N}}V_{\mu_n}$ is ${\cal F}$-acyclic, the covering $U_{\bar{Y}}=\cup_{n\in{\mathbb N}}U_{\mu_n,\rho'}$ is $j^{\dagger}{\cal F}'$-acyclic. It follows $H^q(X,{\cal F}_X)=0=H^q(]\bar{Y}[_{\cal X},j^{\dagger}{\cal F}')=0$ for all $q\ge 2$. 
Also $H^0(X,{\cal F}_X)=H^0(]\bar{Y}[_{\cal X},j^{\dagger}{\cal F}')$ is evident, it remains to show that$$H^1(X,{\cal F}_X)\stackrel{\alpha}{\to}H^1(]\bar{Y}[_{\cal X},j^{\dagger}{\cal F}')$$is bijective. 
This is done using cocycles with respect to the given acyclic coverings: For a cocycle $(g_n)_{n\in{\mathbb N}}\in({\cal F}(V_{\mu_n}))_{n\in{\mathbb N}}$ every $g_n$ extends to some $g_n'\in{\cal F}'(U_{\mu_n,\rho})$ for some $\rho=\rho(g_n)\in\Gamma^*$ with $1<\rho\le\rho'$, that is $g_n'\in j^{\dagger}{\cal F}'(U_{\mu_n,\rho'})$. 
Then $\alpha([(g_n)_n])=[(g_n')_n]$ on the level of cohomologyclasses. 
Now suppose $\alpha([(g_n)_n])=[(g_n')_n]=0$, i.e. there is a cocycle $(f_n)_n\in(j^{\dagger}{\cal F}'(U_{\mu_n,\rho'}))_n$ such that $f_{n+1}|_{U_{\mu_n,\rho'}}-f_n=g_n'$ for all $n$. 
Every $f_n$ is given by an element $f_n\in{\cal F}'(U_{\mu_n,\rho})$ for some $\rho=\rho(f_n)\in\Gamma^*$ with $1<\rho\le\rho'$, therefore induces an element $h_n\in{\cal F}(V_{\mu_{n-1}})$. 
Set $\bar{h}_n=h_{n+1}$ and $\bar{g}_n=g_{n+1}$. The cocycle $(\bar{f}_n)_n$ bounds $(\bar{g}_n)_n$. 
Since $(\bar{g}_n)_n$ and $(g_n)_n$ are cohomologic, it follows that $\alpha$ is injective. Surjectivity: Let $(g_n')_n\in(j^{\dagger}{\cal F}'(U_{\mu_n,\rho'}))_n$ be given. Setting $\bar{g}_n'=g_{n+1}'$, the cocycle $(\bar{g}_n')_n$ is cohomologic to $(g_n')_n$. But every $\bar{g}_n'$ is induced by some $g_{n+1}'\in{\cal F}'(U_{\mu_{n+1},\rho})$ for some $\rho(g_{n+1}')$ with $1<\rho(g_{n+1}')\le\rho'$, inducing also an element $g_n\in{\cal F}(V_{\mu_n})$. One has $\alpha([(g_n)_n])=[(\bar{g}_n')_n]$.\\

\addtocounter{satz}{1}{\bf \arabic{section}.\arabic{satz}}  Similarly, it is not hard to give an interpretation of rigid cohomology with compact support $H_{rig,c}^*(Y/k)$ in terms of de Rham cohomology of dagger spaces; and doing this, we get from \ref{dagpoin} Poincar\'{e} duality for the rigid cohomology of a smooth $Y$, which in \cite{berpoin} is proven in a completely different way, namely by reducing it to Poincar\'{e} duality of cristalline cohomology of smooth proper $\bar{k}$-schemes.\\


\begin{flushleft}
\textsc{Mathematisches Institut der Universit\"at M\"unster\\ Einsteinstrasse 62, 48149 M\"unster, Germany}\\
\textit{E-mail address}: klonne@math.uni-muenster.de
\end{flushleft}
\end{document}